\documentclass[10pt]{amsart} 
\usepackage{graphicx,latexsym,bm,amsmath,amssymb,verbatim,multicol,lscape}
\theoremstyle{definition}

\theoremstyle{remark}

\numberwithin{equation}{section}

\begin{document}
\title{Squares in $(2^2-1)...(n^2-1)$ and $p$-adic valuation}%
\author{Shaofang Hong}
\address{Mathematical College, Sichuan University, Chengdu 610064, P.R. China}
\email{sfhong@scu.edu.cn; sfhong09@gmail.com; s-f.hong@tom.com;
hongsf02@yahoo.com}
\author{Xingjiang Liu}
\address{Mathematical College, Sichuan University, Chengdu 610064, P.R. China}%
\email{xjliuking@163.com}
\thanks{The research of Hong was supported partially by National Science Foundation of China
Grant \# 10971145 and by Program for New Century Excellent Talents
in University Grant \# NCET-06-0785}
\date{\today}%
\begin{abstract}
In this paper, we determine all the squares in the sequence
$\{\prod_{k=2}^n(k^2-1)\}_{n=2}^\infty $. From this, one deduces
that there are infinitely many squares in this sequence. We also
give a formula for the $p$-adic valuation of the terms in this
sequence.
\end{abstract}

\maketitle

\section{Introduction}

The sequence containing infinitely many squares is a common and
interesting subject in number theory.
Cilleruelo [1] proved that there is only one square in the sequence
${\{b_n}\}_{n=1}^{\infty}$, where
$$
b_n:=\prod_{k=1}^n(k^2+1).
$$
In other words, Cilleruelo [1] showed that $b_n$ is a square if and
only if $n=3$. Cilleruelo's proof of this result is analytic. It is
natural to ask the question of determining all the squares in the
analogous sequence ${\{a_n}\}_{n=2}^{\infty}$ defined by
$$a_n:=\prod_{k=2}^n(k^2-1).$$

In this paper, we consider the above question. Namely, we study the
square problem for the sequence ${\{a_n}\}_{n=2}^{\infty}$. Using
the the structure theorem on the solutions of some Pell equations,
we determine all the squares in this sequence. From it one can
conclude that there are infinitely many squares in this sequence. On
the other hand, divisibility properties of integer sequences have
long been objects of interest. $p$-Adic valuation is the modern
language of divisibility. In the last section, we give a formula for
the $p$-adic valuation of $a_n$.

As usual, we let $v_{p}$ be the normalized $p$-adic valuation of
${\bf Q}$. That is, for any positive integer $x$, we have
$v_p(x):=e$ if $p^{e}\parallel x$. We define $s_p(x)$ to be the sum
of the digits of $x$ in its $p$-adic expansion. Namely, we have
$s_p(x):=\sum_{i=0}^l x_i$ if $x=\sum_{i=0}^lx_ip^i$ with $0\le
x_i\le p-1$ for all $0\le i\le l$.

\section{Squares}

In the present section, we determine all the squares in the sequence
$\{a_n\}_{n=2}^\infty $. Let's begin with the following lemma.\\
\\
{\bf Lemma 2.1.} {\it Let $n\geqslant 2.$ Then $a_n$ is a square if
and only if $2n(n+1)$ is a square.}\\
\\
{\it Proof.} Since
$$
a_n=\left(\prod_{k=2}^n(k+1)\right)\left(
\prod_{k=2}^n(k-1)\right)=2n(n+1)\left(\prod_{k=3}^{n-1}k\right)^2,
$$
Lemma 2.1 follows immediately. \hfill$\Box$\\

We need the following well-known result.\\

{\bf Lemma 2.2.} (\cite{[HW]}) {\it All integral solutions of the
Diophantine equation
$$
x^2-2y^2=1
$$
are given by $x+y\sqrt{2}=\pm (1+\sqrt{2})^{2m}$, and all of
$$
x^2-2y^2=-1
$$
by  $x+y\sqrt{2}=\pm (1+\sqrt{2})^{2m+1}$, with $m$ an integer.}\\

We can now give the first main result of this paper.\\
\\
{\bf Theorem 2.3.} {\it Let $n\geqslant 2$. Then $a_n$ is a square
if and only if we have either
$${n=[{k\choose 0}+2{k\choose 2}+...+{2^{(k-1)/2}{k\choose
k-1}}]^2}
$$
with $k\ge 1$ an odd number, or
$${n=2[{k\choose
1}+2{k\choose 3}+...+{2^{(k-2)/2}{k\choose k-1}}]^2}
$$
with $k\ge 2$ an even number.}\\
\\
{\it Proof.} First we show the necessity. Let $a_n$ be a square.
Then by Lemma 2.1, we know that $2n(n+1)$ is a square. So both $2n$
and $n+1$ (resp. both $n$ and $2(n+1)$) are squares if $2|n$ (resp.
$2\not| n$). Let
$$
(a, b):=\left\{\begin{array}{cc}
&(\sqrt{n+1}, \sqrt{\frac{n}{2}}), \ {\rm if} \ 2|n,\\
&(\sqrt{n}, \sqrt{\frac{n+1}{2}}), \ {\rm if} \ 2\not|n.
\end{array}\right.
$$
Then we have $a^2-2b^2=1$ if $2|n$ and $a^2-2b^2=-1$ if $2\not|n$.
Hence by Lemma 2.2, $a+b\sqrt{2}$ has the form $\pm (1+\sqrt{2})^k$,
where $k\in\mathbb{Z}$ is even (resp. odd) if $n$ is even (resp.
odd). But $a$ and $b$ are positive integers. Therefore we must have
$$
a+b\sqrt{2}=(1+\sqrt{2})^k   \eqno (1)
$$
with $k\ge 1$ even (resp. odd) if $n$ is even (resp. odd).

If $n$ is odd, then $k$ is odd and $a=\sqrt{n}$. So $n=a^2$.
Expanding the right hand side of (1) gives us
$$a={k\choose 0}+2{k\choose 2}+...+{2^{(k-1)/2}{k\choose k-1}}.$$
Therefore
$${n=a^2=[{k\choose 0}+2{k\choose 2}+...+{2^{(k-1)/2}{k\choose
k-1}}]^2}.$$

If $n$ is even, then $k$ is even and $b=\sqrt{\frac{n}{2}}$. We then
deduce that $n=2b^2$. Since $k$ is an even number, expanding the
right hand side of (1) we obtain that
$$
b={k\choose
1}+2{k\choose 3}+...+{2^{(k-2)/2}{k\choose k-1}}.
$$
Thus
$$
{n=2b^2=2[{k\choose 1}+2{k\choose 3}+...+{2^{(k-2)/2}{k\choose
k-1}}]^2}
$$
as desired. The necessity is proved.

Conversely, if ${n=[{k\choose 0}+2{k\choose
2}+...+{2^{(k-1)/2}{k\choose k-1}}]^2}$, where $k$ is an odd number,
then $n=[\frac{(1+\sqrt{2})^k+(1-\sqrt{2})^k}{2}]^2$. Therefore
$$\frac{n+1}{2}=\frac{1}{2}[\frac{(1+\sqrt{2})^k-(1-\sqrt{2})^k}{2}]^2. \eqno (2)$$
We expand the right hand side of (2) and then get that
$$\frac{n+1}{2}=[{k\choose
1}+2{k\choose 3}+...+{2^{(k-1)/2}{k\choose k}}]^2,$$ which is a
square. Thus $2n(n+1)$ is a square as required.

If ${n=2[{k\choose 1}+2{k\choose 3}+...+{2^{(k-2)/2}{k\choose
k-1}}]^2}$, where $k$ is an even integer, then in the similar way as
above, we can show that $n+1=[{k\choose 0}+2{k\choose
2}+...+{2^{k/2}{k\choose k}}]^2$. This is a square. Hence $2n(n+1)$
is also a square. Then by Lemma 2.1, the sufficiency follows.

The proof of Theorem 2.3 is complete. \hfill$\Box$\\

From Theorem 2.3 we then deduces that there are infinitely many
squares in the sequence $\{a_n\}_{n=2}^{\infty}$. \\
\\
{\bf Remark.} It would be of interest to consider the square problem
for other similar sequences of positive integers. For example, one
can consider the squares in the sequences
$\{\prod_{k=a+1}^n(k^2-a^2)\}_{n=a+1}^{\infty }$ and
$\{\prod_{k=1}^n(k^2+a)\}_{n=1}^{\infty }$ for any given integer
$a\ge 2$.

\section{The $p$-adic valuation of $a_n$}

In this section, we investigate $p$-adic valuation of $a_n$. First
we give the second main result of this paper.\\
\\
{\bf Theorem 3.1.} {\it Let $n$ be a positive integer. Then each of
the following is true:}

(i).
$$v_2(a_n)=\left\{\begin{array}{cc}
&2n-2-2s_2(\frac{n-1}{2})+v_2(\frac{n+1}{2}), \ {\it if \ n \ is \ odd;}\\
&2n-4-2s_2(\frac{n}{2}-1)+v_2(\frac{n}{2}), \ {\it if \ n \ is \
even.}
\end{array}\right.
$$

(ii). {\it For any odd prime $p$,}

$$v_p(a_n)=v_p(n)+v_p(n+1)+\frac{2}{p-1}(n-1-s_p(n-1)).$$
\\
{\it Proof.} (i). If $n$ is odd, then
$$\begin{array}{rl}
v_2(a_n)&=v_2((3^2-1)(5^2-1)\cdot\cdot\cdot(n^2-1))\\
&=1+2v_2(4\cdot6\cdot8\cdot\cdot\cdot\cdot(n-1))+v_2(n+1)\\
&=1+2(\frac{n-3}{2}+v_2((\frac{n-1}{2})!))+v_2(n+1)\\
&=2n-2-2s_2(\frac{n-1}{2})+v_2(\frac{n+1}{2})
\end{array} \eqno(3)$$
since $k^2-1$ is odd if $k$ is even and Legendre's formula tells us
that $v_2(m!)=m-s_2(m)$ for any positive integer $m$.

If $n$ is even, then $v_2(a_n)=v_2(a_{n-1})$. Hence replacing $n$
with $n-1$ in (3), we obtain
$v_2(a_n)=2n-4-2s_2\left(\frac{n}{2}-1\right)+v_2\left(\frac{n}{2}\right).$
So part (i) is proved.

(ii). Let $p$ be an odd prime. By the Legendre's formula, we get
$$v_p(a_n)=v_p(n)+v_p(n+1)+2v_p((n-1)!)=v_p(n)+v_p(n+1)+\frac{2}{p-1}(n-1-s_p(n-1))$$
as desired. Then part (ii) is proved. The proof of Theorem 3.1 is complete. \hfill$\Box$ \\

Finally, we study the asymptotic behavior of $p$-adic valuation of $a_n$.\\
\\
{\bf Corollary 3.2.} {\it Let $p$ be a prime. Then $v_p(a_n)\sim \frac{2n}{p-1}$ as $n\rightarrow\infty $}.\\
\\
{\it Proof.} By Theorem 3.1 we obtain that
$${\frac{v_2(a_n)}{2n}}=\left\{\begin{array}{cc}
&1-\frac{1}{n}-\frac{s_2(\frac{n-1}{2})}{n}+\frac{v_2(\frac{n+1}{2})}{2n}, \ {\rm if} \ n \ {\rm is \ odd;}\\
&
1-\frac{2}{n}-\frac{s_2(\frac{n}{2}-1)}{n}+\frac{v_2(\frac{n}{2})}{2n},
\ {\rm if} \ n \ {\rm is \ even,}
\end{array}\right.\eqno(4)
$$
and for any odd prime $p$,
$$\frac{v_p(a_n)}{2n/(p-1)}=1-\frac{1}{n}-\frac{s_p(n-1)}{n}
+\frac{p-1}{2}\left(\frac{v_p(n)}{n}+\frac{v_p(n+1)}{n}\right).\eqno(5)$$
Since $v_p(n)\le \log_pn$ and $s_p(n)\le (p-1)(1+\log_pn)$,
Corollary 3.2 follows immediately from (4) and (5). \hfill$\Box$\\

\noindent{\sc Acknowledgement.} The authors would like to thank the
referee for helpful comments.


\begin{thebibliography}{20}

\bibitem {[C]} J. Cilleruelo, {\it Squares in $(1^2+1)\cdot\cdot\cdot(n^2+1)$}, {J. Number
Theory} {\bf 128} (2008), 2488-2491.

\bibitem {[HW]} G.H. Hardy and E.M. Wright, {\it An introduction to the theory of
numbers}, Oxford University Press, 1979.




\end{thebibliography}
\end{document}